\documentclass[a4paper,11pt]{amsart}
\usepackage{a4wide}
\usepackage{amsfonts}
\usepackage{amsmath, amsthm}
\usepackage{amssymb,bbm}

\newcommand{\bF}{\mathbb{F}}

\newcommand{\C}{\mathbb{C}}
\newcommand{\id}{\mathbbm{1}}
\newcommand{\N}{\mathbb{N}}

\newcommand{\Z}{\mathbb{Z}}

\newcommand{\La}{\mathfrak{a}}

\newcommand{\Lf}{\mathfrak{f}}
\newcommand{\Lg}{\mathfrak{g}}
\newcommand{\Lgl}{\mathfrak{gl}}

\newcommand{\Lr}{\mathfrak{r}}
\newcommand{\Ls}{\mathfrak{s}}

\newcommand{\Lz}{\mathfrak{z}}

\newcommand{\rad}{\operatorname{rad}}

\newcommand{\Der}{\operatorname{Der}}
\newcommand{\gDer}{\operatorname{\mathcal{L}Der}}
\newcommand{\gAut}{\operatorname{\mathcal{L}Aut}}

\newcommand{\map}[3]{ #1 : #2 \longrightarrow #3 }
\newcommand{\mapl}[5]{ #1 : #2 \longrightarrow #3 : #4 \longmapsto #5 }
\newcommand{\vb}{\newline \phantom{0} \hfill $\bigtriangleup$}

\newenvironment{bewijs}{ \begin{quote}\begin{proof}[\textsc{Proof}:] }{ \\ \end{proof}\end{quote}}

\theoremstyle{plain} \newtheorem{corollary*}{Corollary}
\theoremstyle{plain} \newtheorem{definition*}{Definition}
\theoremstyle{plain} \newtheorem{lemma*}{Lemma}
\theoremstyle{plain} \newtheorem{proposition*}{Proposition}
\newenvironment{example*}{ \textbf{\textsc{Example. }}  }{  \\ }
\newenvironment{interpretation*}{ \textbf{\text{Interpretation.}}  }{  \\}
\newenvironment{maintheorem*}{ \begin{verse}\textbf{\textsc{Main Theorem}}  }{  \\ \end{verse}}
\newenvironment{proof*}{ \begin{quote}\begin{proof}[\textsc{Proof}:] }{ \\ \end{proof}\end{quote}}
\newenvironment{proposition**}{ \textbf{\textsc{Proposition}} }{ \\ }
\newenvironment{remark*}{ \textbf{\textsc{Remark}}  }{  }
\newenvironment{theorem*}{ \begin{verse}\textbf{\textsc{Theorem}}  }{  \\ \end{verse}}

\title[A characterisation of nilpotent Lie algebras]{A characterisation of nilpotent Lie algebras\\by invertible Leibniz-derivations}
\author{Wolfgang Alexander Moens}
\address{Fakult\"at f\"ur Mathematik\\
Universit\"at Wien\\
  Nordbergstr. 15\\
  1090 Wien \\
  Austria} 
\email {wolfgang.moens@univie.ac.at}
\date{\today}
\thanks{Supported by the FWF grant P21683.}
\keywords{Lie algebra, Leibniz algebra, derivation, pre-derivation, nilpotency, invertibility.}

\subjclass[2000]{17B40, 17B30}

\begin{document}

\begin{abstract} Jacobson proved that if a Lie algebra admits an invertible derivation, it must be nilpotent. He also suspected, though incorrectly, that the converse might be true: that every nilpotent Lie algebra has an invertible derivation. We prove that a Lie algebra is nilpotent if and only if it admits an invertible Leibniz-derivation. The proofs are elementary in nature and are based on well-known techniques.\end{abstract}

\maketitle

\section{Nilpotent Lie algebras and invertible derivations} Which information about the structure of a Lie algebra is contained in its group of automorphisms or in its Lie algebra of derivations? Several partial answers to this question, in terms of sufficient conditions for the nilpotency or solvability of the Lie algebra, are known - see for example the theorems by Borel, Serre and Mostow \cite{BorelSerre,BorelMostow}. We will mainly be interested in the following well-known theorem by Jacobson on the invertibility of Lie algebra derivations.

\begin{theorem*}[\cite{Jacobson}] A Lie algebra over a field of characteristic zero is nilpotent if it admits an invertible derivation. \end{theorem*}

In that same paper, Jacobson also inquired whether the converse statement were true. Dixmier and Lister demonstrated in \cite{DixmierLister} that it was not the case by presenting nilpotent Lie algebras admitting only nilpotent derivations. They called such Lie algebras \emph{characteristically nilpotent}. For an overview of the study of characteristically nilpotent Lie algebras, see \cite{AncocheaCampoamor}. In his study of metrics on Lie groups \cite{Mueller}, M\"uller introduced a natural generalisation of derivations: the so-called pre-derivations. These generalised derivations allowed Bajo to take Jacobson's theorem one step further (cf. proposition $2.8$ in \cite{Burde}).

\begin{theorem*}[\cite{Bajo}] A Lie algebra over a field of characteristic zero is nilpotent if it admits an invertible pre-derivation. \end{theorem*}

Once again, it was asked whether the converse were true; but counterexamples were found: nilpotent Lie algebras admitting only nilpotent pre-derivations. Such Lie algebras are called \emph{strongly nilpotent}, \cite{Burde}. By stating Jacobson's theorem even more generally, in terms of Leibniz-derivations, the result still holds and its converse becomes true. The main aim of this note is to prove the following theorem.

\begin{maintheorem*} A Lie algebra over a field of characteristic zero is nilpotent if and only if it has an invertible Leibniz-derivation. \end{maintheorem*}

This note is organised as follows. We first introduce the algebras of Leibniz-derivations. We then present some existence results under the assumption that the Lie algebra is nilpotent. After this, we give some sufficient conditions under which certain derivation algebras coincide and we consider a stronger notion of invariance of ideals. Finally, we illustrate how the results in the preceding paragraphs can be used to generalise Bajo's proofs. Unless it is explicitly stated otherwise, every vector space is assumed to be finite-dimensional over an algebraically closed field $\bF$ of characteristic zero.

\section{Leibniz Algebras and Leibniz derivations} In this paragraph we consider a straightforward generalisation of Leibniz-algebras: the $k$-Leibniz algebras as introduced in \cite{CasasLodayPirashvili}. Let us briefly recall some relevant definitions and examples in this context, but note that we use a slightly different parametrisation than the one in \cite{CasasLodayPirashvili}. Consider a vector space $V$ and a natural number $k \geq 1$. A $k$-Leibniz algebra structure on $V$ is a $(k+1)$-linear map $V \times \cdots \times V \longrightarrow V : (x_1,\ldots,x_{k+1}) \longmapsto [x_1,\ldots,x_{k+1}]$ satisfying the identity, $$ [[x_1,\ldots,x_{k+1}],y_1,\ldots,y_k] = \sum_{i=1}^n [x_1,\ldots,x_{i-1},[x_i,y_1,\ldots,y_k],x_{i+1},\ldots,x_{k+1}] , $$ for all $x_1,\ldots,x_{k+1};y_1,\ldots,y_k$ in $V$. Lie algebras, and more generally Leibniz-algebras, are examples of $1$-Leibniz algebras. Lie triple systems are examples of $2$-Leibniz algebras. If $\Lg$ is a Leibniz algebra, $k$ is a natural number at least one, and $x_1,\ldots,x_{k+1}$ are elements of $\Lg$, we may define the multilinear bracket $[x_1,x_2,\ldots,x_{k+1}]_{k+1}$, or when the length of the bracket is clear from the context, $[x_1,x_2,\ldots,x_{k+1}]$, as the nested expression $$[x_1,[x_2,[x_3,\ldots ,[x_{k-1},[x_{k},x_{k+1}]] \ldots ]]].$$ This defines a $k$-Leibniz-algebra structure structure on $\Lg$, denoted by $\mathcal{L}_{k}(\Lg)$. Many more examples can be found in the literature. Subalgebras, ideals and solvability of Leibniz-algebras are defined in the natural way. A subspace $\mathcal{S}$ of a Leibniz $k$-algebra $\mathcal{L}$ is a \emph{subalgebra} of $\mathcal{L}$, if $[\mathcal{S},\ldots,\mathcal{S}]_{k+1} \subseteq \mathcal{S}$. A subspace $\mathcal{I}$ is an \emph{ideal}, if it satisfies $ [\mathcal{I},\mathcal{L},\ldots,\mathcal{L}]_{k+1} + \ldots + [\mathcal{L},\ldots,\mathcal{L},\mathcal{I}]_{k+1} \subseteq \mathcal{I} $. The Leibniz algebra $\mathcal{L}$ is \emph{solvable}, if the derived sequence,
$$  \mathcal{S}_{1}(\mathcal{L}) = \mathcal{L}, \; 
 \mathcal{S}_{t+1}(\mathcal{L}) = [\mathcal{S}_{t}(\mathcal{L}),\ldots,\mathcal{S}_{t
 }(\mathcal{L})]_{k+1}, \text{ for all $t \geq 1$}, $$
descends to $\{0\}$.
We now introduce a generalisation of Lie algebra (pre)-derivations. It is important to note that many different generalisations of algebra derivations can be found in the literature, each emphasizing a particular property of classical derivations.

\begin{definition*} Let $k \in \N_0$. A \emph{Leibniz-derivation} of \emph{order} $k$ for a Lie algebra $\Lg$ is an endomorphism $P$ of that Lie algebra satisfying the identity \begin{eqnarray*} P([x_1, \ldots , x_{k+1}]) &=& [P(x_1),x_2,\ldots,x_{k+1}] + [x_1,P(x_2),\ldots,x_{k+1}] \\ & & + \ldots + [x_1,x_2,\ldots,P(x_{k+1})], \end{eqnarray*} for all $x_1,\ldots,x_{k+1}$ in $\Lg$. Let $\gDer_k(\Lg)$ be the set of all Leibniz-derivations for $\Lg$ of order $k$ and let $\gDer(\Lg)$ be the set of all Leibniz-derivations: $\bigcup_{k \in \N_0} \gDer_k(\Lg) $. \end{definition*}

\begin{interpretation*} The derivations of $\Lg$ are exactly the Leibniz-derivations of order $1$. The pre-derivations are those of order $2$. In \cite{JacobsonTriple}, Jacobson introduced Lie triple systems. The pre-derivations of $\Lg$ can also be interpreted as the derivations of the Lie triple system induced by $\Lg$. Even more generally, we may consider the $k$-Leibniz algebra $\mathcal{L}_{k}(\Lg)$ that is naturally induced by $\Lg$. The derivations of $\mathcal{L}_{k}(\Lg)$ correspond with the Leibniz-derivations for $\Lg$ of order $k$: $$ \Der(\mathcal{L}_{k}(\Lg)) = \gDer_{k}(\Lg) .$$ This suggests that Bajo's theorem is a natural generalisation of Jacobson's, and that our main theorem is a natural generalisation of Bajo's. Note that the order of a Leibniz-derivation is not unique:\end{interpretation*}

\begin{lemma*} If $s,t \in \N_0$ and $s | t$, then $\gDer_s(\Lg) \subseteq \gDer_t(\Lg).$ If $k,l \in \N_0$, then $$\gDer_k(\Lg) \cap \gDer_l(\Lg) \subseteq \gDer_{k+l}(\Lg).$$  \end{lemma*}

\begin{proof*} $(i.)$ Suppose $P \in \gDer_s(\Lg)$. By recursively applying the definition (of a Leibniz-derivation of order $s$) $\frac{t}{s}$ times to a bracket of the form $[x_1,\ldots,x_{t+1}]$, we see that also $P \in \gDer_t(\Lg)$. $(ii.)$ Similarly, we may write a bracket of the form $[x_1,\ldots,x_{k+l+1}]$ as $[x_1,\ldots,x_k,[x_{k+1},\ldots,x_{k+(l+1)}]]$. \end{proof*}

\begin{proposition*} The subsets $\gDer(\Lg)$ and all $\gDer_k(\Lg)$ of $\Lgl(\Lg)$ are in fact subalgebras and we have the chain $$\operatorname{Inn}(\Lg) \subseteq \Der(\Lg) \subseteq \gDer_k(\Lg) \subseteq \gDer(\Lg) \subseteq \Lgl(\Lg).$$\end{proposition*}

\begin{proof*} It can be checked that each $\gDer_k(\Lg)$ is a vector space. If $P \in \gDer_k(\Lg)$ and $Q \in \gDer_l(\Lg)$ then $P,Q \in \gDer_{(k,l)}(\Lg)$ by the previous lemma. This shows that also $\gDer(\Lg)$ is a vector space. Because of the previous lemma, it now suffices to prove that $[\gDer_m(\Lg),\gDer_m(\Lg)] \subseteq \gDer_m(\Lg)$. Let $P,Q \in \gDer_m(\Lg)$ and $x_1,\ldots,x_{m+1} \in \Lg$. Then $[P,Q]([x_1,\ldots,x_{m+1}])$ is equal to 
\begin{eqnarray*}
\sum_{i \neq j} [x_1,\ldots,P(x_i),\ldots,Q(x_j)
 \ldots,x_{m+1}]
 &+& \sum_l [x_1,\ldots,(P \circ Q)(x_l),\ldots,x_{m+1}] \\
 - \sum_{i \neq j} [x_1,\ldots,Q(x_i),\ldots,P(x_j),\ldots,x_{m+1}]
 &-& \sum_l [x_1
 ,\ldots,(Q \circ P)(x_l),\ldots,x_{m+1}],
\end{eqnarray*} which in turn is equal to $\sum_l [x_1,\ldots,[P,Q](x_l),\ldots,x_{m+1}]$. This finishes the proof. \end{proof*} $ $

\begin{remark*} This chain collapses for Lie algebras with special properties. For nilpotent Lie algebras, the first three terms are generically distinct and the last two terms coincide (section 3). For perfect Lie algebras without center and for reductive Lie algebras, the middle three terms coincide (section $4$). \end{remark*} \newline

The Lie algebras $\gDer_k$ and $\gDer$ are clearly invariants of Lie algebras. We now apply a well-known technique concerning the behaviour of (Lie) algebra invariants under degenerations. We use the following fact: if a sequence of matrices of nullity at least $n_0$ converges (in the usual topology), then the nullity of the limit is also at least $n_0$.

\begin{proposition*} If $\Lg \rightarrow \Lg_0$ is a degeneration of \emph{complex} Lie algebras and $k \in \N_0$, then $$\dim(\gDer_k(\Lg)) \leq \dim(\gDer_k(\Lg_0)).$$\end{proposition*}

\begin{proof*} Fix the natural number $k \geq 1$. Choose a basis for the underlying vector space, say $V$, and let $n$ be its dimension. With respect to this basis, we may identify the endomorphisms of $V$ with the $(n \times n)$-matrices and the Lie algebra structures on $V$ with their structure constants. One can then explicitly construct (cf. \cite{BenesBurde,NesterenkoPopovych}) a polynomial map $\map{M}{ \mathcal{L}(V) \subseteq \C^3 }{\operatorname{Mat}({n^3} , {n^2})}$ from the Lie algebra laws on $V$ to the $(n^2 \times n^3)$-matrices satisfying the following property: if $\mu \in \mathcal{L}(V)$, then $\gDer_k(\mu) $ and $ \ker(M(\mu))$ are isomorphic as vector spaces. This implies that $\dim(\gDer_k(\mu))$ is equal to the nullity of $ M(\mu) $.\newline \\* Now suppose that $\mu \in \mathcal{L}(V)$ degenerates to $\mu_0 \in \mathcal{L}(V)$. Then (cf. \cite{BenesBurde,NesterenkoPopovych}) $\mu \in \mathcal{L}(V)$ sequentially contracts to $\mu_0 \in \mathcal{L}(V)$, say $\lim_{\varepsilon \rightarrow 0} ( \mu_\varepsilon ) = \mu_0.$ By the above remarks, it suffices to prove $(i.)$ that each $M(\mu_\varepsilon)$ has the same nullity as $M(\mu)$ and $(ii.)$ that  $\lim_{\varepsilon \rightarrow 0} (M(\mu_\varepsilon)) = M(\mu_0)$. The former follows from the fact that $\mu$ is, by definition, isomorphic to each $\mu_\varepsilon$. Since the map $M$ is polynomial, and thus continuous, we also have $\lim_\varepsilon (M(\mu_\varepsilon)) = M(\lim_\varepsilon \mu_\varepsilon) = M(\mu_0).$ This finishes the proof.\end{proof*}

\begin{remark*} This proposition generalises the well-known result for $k = 1$. Note, however, that (unlike in the $(k=1)$-case) a proper degeneration \emph{need not} produce a strict inequality for higher $k$. Counterexamples can already be found in dimension 7.
\end{remark*}

\section{Nilpotent Lie algebras} Now suppose that $\Lg$ is a nilpotent Lie algebra of class $c > 0$. Then it is easy to construct Leibniz-derivations. In fact, \emph{every} endomorphism of $\Lg$ is a Leibniz-derivation of order $c$ (and all higher orders). The converse is just as easy to prove:

\begin{proposition*} Consider a Lie algebra $\Lg$ and a $k$ in $\N_0$. Then $\Lg^k = 0$ if and only if $\gDer_k(\Lg)=\Lgl(\Lg)$. In particular, the minimal $k$ for which the equalities hold, is the nilpotency class of $\Lg$. \end{proposition*}

From this it is clear that invertible endomorphisms of $\Lg$ are Leibniz-derivations of arbitrarily high order. It is not clear however, whether they can also have low orders. We generalise the following existence results \cite{Burde}: every nilpotent Lie algebra of class at most two (four) has an invertible (pre)derivation.

\begin{proposition*} Every nilpotent Lie algebra of nilpotency class $c$ has a Leibniz-derivation of order $\lceil \frac{c}{2} \rceil$.\end{proposition*}

\begin{bewijs} Let $\Lg$ be the Lie algebra and set $q = \lceil \frac{c}{2} \rceil$. Choose a vector subspace $W$ of $\Lg$ complementary to $\Lg^{q+1}$: $\Lg = W + \Lg^{q+1}$. Define the map $P$ by its restrictions to $W$ and $\Lg^{q+1}$: $P|_{W} = \id_W$ and $P|_{\Lg^{q+1}} = (q+1) \id_{\Lg^{q+1}}$. Then it is easy to check that $P$ is a Leibniz-derivation for $\Lg$ of order $q$. Note that it is even semisimple. \end{bewijs}

The following proposition by Bajo, Benayadi and Medina characterises the Lie algebras admitting an invertible derivation, $k=1$. We are not aware of an analogous description for higher values of $k$. \newline

\begin{proposition**}[\cite{BajoBenayadiMedina}] \emph{A Lie algebra over a field of characteristic zero admits an invertible derivation if and only if it is the quotient of a Lie algebra that is quadratic and symplectic by an ideal that is Lagrangian and completely isotropic.} \end{proposition**}

We have already remarked that certain derivation algebras may coincide. We can also consider the converse in the case of nilpotent Lie algebras. A classical result by Dixmier and Schenkman (\cite{Dixmier} and \cite{Jacobson}) implies that the inequality $\operatorname{Inn}(\Lg) \subseteq \Der(\Lg)$ is in fact strict. If the nilpotency class is at least two, then also $\Der(\Lg) \varsubsetneq \operatorname{pDer}(\Lg) = \gDer_2(\Lg)$, \cite{Burde}. This last statement can be generalised:

\begin{proposition*} Consider a nilpotent Lie algebra $\Lg$ of class $c$. Let $k$ and $l$ be natural numbers satisfying $k \equiv 1 \operatorname{mod} l$
%$l | k-1$
and $k \leq c$. Then we can construct an element of $\gDer_k(\Lg) \setminus \gDer_l(\Lg)$. \end{proposition*}

\begin{proof*} Choose elements $e_1,\ldots,e_k$ of $\Lg^1$ such that $u = [e_1,\ldots,e_k] \in \Lg^{k} \setminus \Lg^{k+1}$. Choose a basis $\mathcal{B}$ for $\Lg$ containing $u$ such that $\langle \mathcal{B} \setminus u \rangle$ contains $e_1,\ldots,e_k$ and $\Lg^{k+1}$. Choose a non-zero central element $z$ of $\Lg$. Then define the linear map $\map{P_z}{\Lg}{\Lg}$ on this basis (and extend it linearly): $P_z(u) = z$ and $P_z$ is zero on all other basis vectors. Then it is clear that $P_z$ is a Leibniz-derivation for $\Lg$ of order $k$. But $P_z$ does not have order $l$: otherwise, it would also have order $k-1$ by lemma $1$ and we would obtain the contradiction $$ z = P_z(u) = [P_z(e_1),\ldots,e_{k}] + \ldots + [e_1,\ldots,P_z(e_k)] = 0 .$$\qedhere \end{proof*}

\section{Perfect and reductive Lie algebras} We now show that the Lie algebras $\gDer(\Lg)$ and $\Der(\Lg)$ coincide if $\Lg$ is reductive or perfect. Recall that the Lie algebra is called reductive if its radical coincides with its center, and it is called perfect if the Lie algebra coincides with its commutator.

\begin{lemma*} Consider a perfect Lie algebra $\Lg$. % with trivial center.
Choose a natural number $k \geq 2$ and an element $m \in \bF^\times \setminus \Z^-$. Among all the endomorphisms $\map{f}{\Lg}{\Lg}$ of $\Lg$, only the zero-map satisfies the identity,
\begin{eqnarray*}m \cdot f([x_1,\ldots,x_k]) &+& [x_1,f([x_2,\ldots,x_k])] \\ &+& [x_1,x_2,f([x_3,\ldots,x_k])] \\ &+& \ldots \\ &+& [x_1,x_2,\ldots,x_{k-1},f(x_k)] \\ &=& 0.\phantom{O}(\ast_{m,k})\end{eqnarray*}\end{lemma*}

\begin{proof*} We prove the lemma by induction on the parameter $k$. Consider the substitution $x_1 \mapsto [t,x_1]$ in the identity $(\ast_{m,k})$. We apply the Jacobi-identity to the left part of each term but the first and then twice the identity $(\ast_{m,k})$ to obtain $$ - f([[t,x_1],x_2,\ldots,x_k]) + [f([t,x_2,\ldots,x_k]),x_1] + [t,f([x_1,x_2,\ldots,x_k])] = 0. $$ Now group $x_{k-1}$ and $x_k$ together so that we may once again apply identity $(\ast_{m,k})$:
\begin{eqnarray*}
&-& f([[t,x_1],x_2,\ldots,x_k]) \\
&+& m \cdot f([x_1,t,x_2,\ldots,x_k]) - m \cdot f([t,x_1,x_2,\ldots,x_k]) \\
&+& [x_1,t,f([x_2,\ldots,x_k])] - [t,x_1,f([x_2,\ldots,x_k])] \\
&+& [x_1,t,x_2,f([x_3,\ldots,x_k])] - [t,x_1,x_2,f([x_3,\ldots,x_k])]\\
&+& \cdots \\
&+& [x_1,t,x_2,\ldots,x_{k-2},f([x_{k-1},x_k])] - [t,x_1,x_2,\ldots,x_{k-2},f([x_{k-1},x_k])] \\
&=& 0.
\end{eqnarray*}
Note that $[[t,x_1],p] = [t,x_1,p] - [x_1,t,p]$ for all $p$ in $\Lg$, so that the first three terms add up to $(m+1) \cdot f([[x_1,t],x_2,\ldots,x_k])$. The remaining terms can also be grouped together and we obtain
\begin{eqnarray*}
 (m+1) \cdot f([[x_1,t],x_2,\ldots,x_k]) &+& [[x_1,t],f([x_2,\ldots,x_k])] \\
&+& [[x_1,t],x_2,f([x_3,\ldots,x_k])] \\
&+& \cdots \\
&+& [[x_1,t],x_2,\ldots,x_{k-2},f([x_{k-1},x_k])] \\
&=& 0.
\end{eqnarray*}
Note that $m+1$ is non-zero. If $k = 2$, this identity reduces to $(m + 1) \cdot f([[x_1,t],x_2]) = 0$ and we may conclude that $f$ vanishes. This proves the base of the induction. Now assume that $k \geq 3$. Since $\Lg$ is generated by its commutators, we may perform the substitutions $[t,x_1] \mapsto x_1$ and $[x_{k-1},x_k] \mapsto x_{k-1}$.
So we observe that $f$ also satisfies the identity $(\ast_{m+1,k-1})$. Since $m+1$ belongs to $\bF^\times \setminus \Z^-$ and since $k-1 \geq 2$, we may apply the induction hypothesis, which guarantees the vanishing of $f$. This finishes the induction. \end{proof*}

\begin{proposition*} All Leibniz-derivations of a perfect Lie algebra $\Lg$ are derivations in the classical sense: $\gDer(\Lg) = \Der(\Lg)$. \end{proposition*}

\begin{proof*} Let $P$ be a Leibniz-derivation of the perfect and centerless Lie algebra $\Lg$, touching $\Lg$ at $k \geq 1$. So it suffices to prove that the alternating, bilinear map $$\mapl{\omega}{\Lg \times \Lg}{\Lg}{(x;y)}{P([x,y]) - [P(x),y] - [x,P(y)]}$$ vanishes. This is true - by definition - for $k=1$. So we may assume that $k \geq 2$. We can rewrite the definition of a Leibniz-derivation in terms of this $\omega$: for all $x_1,\ldots,x_{k+1}$ in $\Lg$, we have
\begin{eqnarray*}
\omega(x_1;[x_2,\ldots,x_{k+1}]) &+& [x_1,\omega(x_2;[x_3,\ldots,x_{k+1}])] \\ &+& [x_1,x_2,\omega(x_3;[x_4,\ldots,x_{k+1}])] \\ &+& \cdots \\ &+& [x_1,x_2,\ldots,x_{k-1},\omega(x_k;x_{k+1})] \\ &=& 0. \phantom{O} (\dagger) \label{Recursion}
\end{eqnarray*}
Note that we cannot apply the technical lemma to this identity ($\omega$ is defined on $\Lg \times \Lg$). But we can explicitly construct a linear map $\map{f}{\Lg}{\Lg}$ satisfying the identity $ \omega(x;y) = [x , f(y)]$ for all $x,y$ in $\Lg$. We do this as follows.  Choose an ordered basis $(z_1,\ldots,z_d)$ for $\Lg$. Since $\Lg$ is perfect, each of these basis vectors $z_i$ can be written as a sum of $k$-commutators: $$z_i = \sum_r [y_1^i(r),\ldots,y_k^i(r)] = \sum_r y(i,r)$$ for some $y_j^i(r)$ in $\Ls$. We should remark that this decomposition is not necessarily unique. Now define $f(z_i)$ as $$ \sum_r \left( P(y(i,r)) - [P(y_1^i(r)),y_2^i(r),\ldots,y_k^i(r)] - \ldots - [y_1^i(r),\ldots,P(y_k^i(r))] \right) ,$$ and extend $f$ linearly to all of $\Lg$. By using the definition of a Leibniz-derivation, one can then check that this function satisfies the identity $\omega(x;y) = [x , f(y)]$. We may then rewrite the identity $(\dagger)$ to obtain
\begin{eqnarray*}
[x_1,f([x_2,\ldots,x_{k+1}])] &+& [x_1,x_2,f([x_3,\ldots,x_{k+1}])] \\ &+& [x_1,x_2,x_3,f([x_4,\ldots,x_{k+1}])] \\ &+& \cdots \\ &+& [x_1,x_2,\ldots,x_{k-1},x_k,f(x_{k+1})] \\ &=& 0. \label{Recursion}
\end{eqnarray*}
Note that $f$ leaves the center of $\Lg$ invariant: $0 = [\Lz(\Lg),f(\Lg)] = \omega(\Lz(\Lg);\Lg) = - \omega(\Lg;\Lz(\Lg)) = - [\Lg,f(\Lz(\Lg))]$. This $f$ then naturally projects to an endomorphism $\overline{f}$ of the perfect Lie algebra $\overline{\Lg} = \Lg \slash \Lz(\Lg)$. So we obtain the identity \begin{eqnarray*}
\overline{f}([\overline{x}_2,\ldots,\overline{x}_{k+1}]) &+& [\overline{x}_2,\overline{f}([\overline{x}_3,\ldots,\overline{x}_{k+1}])] \\ &+& [\overline{x}_2,\overline{x}_3,\overline{f}([\overline{x}_4,\ldots,\overline{x}_{k+1}])] \\ &+& \cdots \\ &+& [\overline{x}_2,\ldots,\overline{x}_{k-1},\overline{x}_k,\overline{f}(\overline{x}_{k+1})] \\ &=& \overline{0}. \label{Recursion}
\end{eqnarray*}
This means that $\overline{f}$ satisfies the identity $(\ast_{m,k})$ of the previous lemma (with $m = 1$ and $k \geq 2$). We may then conclude that $f$ maps into the center of $\Lg$, and that that $\omega$ vanishes: $\omega(\Lg;\Lg) = [\Lg,f(\Lg)] \subseteq [\Lg,\Lz(\Lg)] = 0$. This finishes the proof. \end{proof*} $ $

If the condition fails, the equality $\Der(\Lg) = \gDer(\Lg)$ may still hold. Examples are given by the abelian Lie algebras $\C^d$ and the $2$-dimensional affine Lie algebra $\La \Lf \Lf_1(\C)$.

\begin{corollary*} All Leibniz-derivations of a semisimple Lie algebra are inner and thus singular. \end{corollary*}

The following lemma is well-known for $k=1$ (\cite{Togo}, Lemma $1$) and the proof for higher $k$ should be obvious.

\begin{lemma*} Consider a Lie algebra $\Lg$ and suppose it can be decomposed as a direct sum of ideals: $\Lg = \oplus_i \Lg_i$. Define $ \gDer_k(\Lg_i,\Lg_i) = \gDer_k(\Lg_i) $ and, for $i\neq j$, let $\gDer_k(\Lg_i,\Lg_j)$ consist of the endomorphisms $f$ of $\Lg$ satisfying the following conditions: $ f(\Lg_l) = 0$ if $l \neq i$, $f(\Lg) \subseteq \Lg_j$, $f([\Lg_i,\ldots,\Lg_i]_{k+1}) = 0$ and $[f(\Lg_i),\Lg_j,\ldots,\Lg_j]_{k+1} + \ldots + [\Lg_j,\ldots,\Lg_k,f(\Lg_i)]_{k+1} = 0.$ Then $$\gDer_k(\Lg) = \sum_{i,j} \gDer_k(\Lg_i,\Lg_j).$$ \end{lemma*}

\begin{proposition*} All Leibniz-derivations of a reductive Lie algebra $\Lg$ are derivations in the classical sense: $\gDer(\Lg) = \Der(\Lg).$ \end{proposition*}

\begin{proof*} We may decompose $\Lg$ as $\Ls \oplus \La$, with $\Ls$ semisimple and $\La$ abelian. Since $\Ls$ is perfect and centerless, the terms $\gDer_k(\Ls,\La)$ and $\gDer_k(\La,\Ls)$ vanish. By the previous lemma, we obtain $\gDer_k(\Lg) = \gDer_k(\Ls) + \gDer_k(\La)$. Alternatively, we observe that perfect ideals are invariant under Leibniz-derivations. By proposition 9, also the radical $\La$ is invariant. This again gives us $ \gDer_k(\Lg) = \gDer_k(\Ls) + \gDer_k(\La) $. \newline \\* Observe that $ \gDer_k(\Ls) $ is equal to $\operatorname{Inn}(\Lg)$ by corollary 1. For the abelian Lie algebra $\La$ we have $\gDer_k(\La) = \Lgl(\La)$. We may then combine these observations to conclude that $\gDer(\Lg) = \operatorname{Inn}(\Ls) + \Lgl(\La) = \Der(\Lg).$ \qedhere\end{proof*}

\section{Invariance of the solvable radical} The ideals of a Lie algebra $\Lg$ that are invariant under all classical derivations are called \emph{characteristic} ideals. Obvious examples are the terms of the descending central series, the terms of the derived central series, and the terms of the ascending central series: $\Lg^k$, $\Lg^{(k)}$ resp. $\Lz_k(\Lg)$. The maximal nilpotent ideal $\operatorname{nil}(\Lg)$ and the maximal solvable ideal $\rad(\Lg)$ are also characteristic ideals. In this paragraph, we show that the solvable radical (of $\Lg$) is also invariant under all Leibniz-derivations. \newline

We define the \emph{radical} $\rad(\mathcal{L}_k(\Lg))$ of $\mathcal{L}_k(\Lg)$ to be the (unique) maximal solvable ideal of $\mathcal{L}_k(\Lg)$ Note that this radical is well defined, since the sum of two solvable ideals of $\mathcal{L}_k(\Lg)$ is again a solvable ideal of $\mathcal{L}_k(\Lg)$.

\begin{proposition*} Let $\Lg$ be a Lie algebra and $k \in \N_0$. Then $ \rad(\mathcal{L}_k(\Lg)) = \rad(\Lg).$ \end{proposition*}

\begin{proof*} Any solvable ideal of $\Lg$ will also be a solvable ideal of $\mathcal{L}_k(\Lg)$, so that it suffices to prove the inclusion $\rad(\mathcal{L}_k(\Lg)) \subseteq \rad(\Lg)$. Let $X$ be a solvable ideal of $\mathcal{L}_k(\Lg)$. Consider a Levi-decomposition $\Lg = \Ls \ltimes \Lr$ for $\Lg$ and let $\map{\pi}{\Lg}{\Ls}$ be the natural quotient map. Then this map is a morphism of Leibniz-$k$-algebras so that $\pi(X)$ is solvable as a $k$-Leibniz algebra. But $\pi(X)$ is also an ideal of $\Ls$: $[X,\Lg,\ldots,\Lg]_{k+1} \subseteq X \Rightarrow [\pi(X),\Ls] \subseteq \pi(X)$. This implies that $\pi(X)$ is semi-simple and thus perfect. Only $\pi(X) = 0$ is both solvable (as a $k$-Leibniz algebra) and equal to its own commutator so that we conclude that $X$ is contained in $\rad(\Lg)$. This finishes the proof. \end{proof*}

\begin{proposition*} The radical of a Lie algebra is invariant under all Leibniz-derivations.
\end{proposition*}

\begin{proof*} This proof is a generalisation of one by Hochschild \cite{Hochschild}. Let $\Lg$ be the Lie algebra, $k \in \N_0$. By the previous proposition, we may identify the radicals $\operatorname{rad}(\Lg)$ and $\operatorname{rad}(\mathcal{L}_k(\Lg))$. Denote them by $\Lr$. Consider the derived sequence of $\mathcal{L}_k(\Lr)$: $$ \Lr = \Lr_1 \varsupsetneq \Lr_2 \varsupsetneq \cdots \varsupsetneq \Lr_p \varsupsetneq \Lr_{p+1} = 0, $$ where $\mathcal{S}_t(\mathcal{L}_k(\Lr)) := \Lr_t$. Note that the terms are ideals of $\mathcal{L}_k(\Lr)$ since they are even ideals of $\Lr$. Consider a derivation $\delta$ of $\mathcal{L}_k(\Lg)$. We will then show, that $$\delta^i(\Lr_t) \subseteq \Lr,$$ for all $i \in \N_0$ and all $1 \leq t \leq p+1$ by using induction on $t$. The basis of the induction is given by $t = p+1$. The induction hypothesis is then $ \delta^i(\Lr_{t+1})  \subseteq \Lr $ for all $i$. We now show that also $ \delta^i(\Lr_{t})  \subseteq \Lr $ holds for all $i$. \newline

The set $A = \Lr + \delta(\Lr_t)$ is an ideal of $\mathcal{L}_t(\Lg)$ since
\begin{eqnarray*}
 [\Lg,\ldots,\Lg,A,\Lg,\ldots,\Lg]_{k+1} &\subseteq& [\Lg,\ldots,\Lg,\Lr ,\Lg,\ldots
 ,\Lg]_{k+1} \\ && + [\Lg,\ldots,\Lg,\delta(\Lr_t),\Lg,\ldots,\Lg]_{k+1} \\
  &\subseteq& \Lr + \Lr_t + \delta(\Lr_t) \\
  &=& A.
\end{eqnarray*}
The derived algebra $[A,\ldots,A]_{k+1}$ of $A$ is contained in $\Lr$:
\begin{eqnarray*}
 [A,\ldots,A]_{k+1}
 &\subseteq& \Lr + [\delta(\Lr_t),\ldots,\delta(\Lr_t)]_{k+1} 
 \\
 &\subseteq& \Lr + \Lr + \delta^{k+1}(\Lr_{t+1}) \\
 &\subseteq& \Lr,
\end{eqnarray*}
by the induction hypothesis. We may then conclude that $A$ is a solvable ideal of $\mathcal{L}_t(\Lg)$ and thus contained in $\Lr$. In particular, we have $\delta(\Lr_{t}) \subseteq \Lr$. Note that the proof does not end here. Now suppose that $\delta^i(\Lr_t) \subseteq \Lr$ for all $1 \leq i < n$. We will then show that also $\delta^n(\Lr_t) \subseteq \Lr$. \newline

As before, we can show that the set $B = \Lr + \delta^n(\Lr_t)$ is an ideal of $\mathcal{L}_k(\Lg)$:
\begin{eqnarray*}
[\Lg,\ldots,\Lg,B,\Lg,\ldots,\Lg]_{k+1} &\subseteq& \Lr + [\Lg,\ldots,\Lg,\delta^n(\Lr_t),\Lg,\ldots,\Lg]_{k+1} \\
&\subseteq& \Lr + \Lr_t + \delta^n(\Lr_t) \\
&=& B,
\end{eqnarray*}
by the hypothesis $\delta^i(\Lr_t) \subseteq \Lr$ for all $i < n$. Again, we can show that the derived algebra of $B$ is in $\Lr$:
\begin{eqnarray*}
[B,\ldots,B]_{k+1} &\subseteq& \Lr + [\delta^n(\Lr_t),\ldots,\delta^n(\Lr_t)]_{k+1} \\
&\subseteq& \Lr + \Lr + \delta^{n(k+1)}(\Lr_{t+1}) \\
&\subseteq& \Lr + \Lr = \Lr,
\end{eqnarray*}
by the hypothesis that $\delta^i(\Lr_{t+1}) \subseteq \Lr$ for all $i$. As before, we conclude that $B \subseteq \Lr$ and thus also $\delta^n(\Lr_t) \subseteq \Lr$. This finishes the proof.
\end{proof*}

Though we will not be needing it later, we mention that one can also prove the invariance of the nilradical under all Leibniz-derivations. For this, we consider the automorphisms of a Leibniz-algebra.

\begin{definition*} Let $k \in \N_0$. A \emph{Leibniz} $k$-\emph{automorphism} for a Lie algebra $\Lg$ is an invertible endomorphism $A$ of that Lie algebra satisfying the identity $$ A([x_1, \ldots , x_{k+1}]) = [A(x_1),A(x_2),\ldots,A(x_{k+1})], $$ for all $x_1,\ldots,x_{k+1}$ in $\Lg$. Let $\gAut_k(\Lg)$ be the set of all Leibniz $k$-automorphisms and let $\gAut(\Lg)$ be the set of all Leibniz automorphisms: $\bigcup_{k \in \N_0} \gAut_k(\Lg).$\end{definition*}

\begin{interpretation*} As in the case of derivations, the automorphisms of $\Lg$ correspond with $k = 1$. M\"uller's pre-automorphisms correspond with $k = 2$. For each Lie algebra $\Lg$ and each natural number $k \geq 1$, $\gAut_k(\Lg) = \operatorname{Aut}(\mathcal{L}_{k}(\Lg))$ is a Lie group with Lie algebra $\gDer_k(\Lg)$.\end{interpretation*}

\begin{proposition*} The solvable radical and the nilradical of a Lie algebra are invariant under all Leibniz-automorphisms and all Leibniz-derivations. \end{proposition*}

\begin{proof*} For $k=1$: the result is well-known for Lie algebras. For $k=2$: the result for pre-automorphisms is due to M\"uller (\cite{Mueller}, propositions $3.3$ and $3.4$), the one for pre-derivations is due to Bajo (\cite{Bajo}, lemma $1$). For $k\geq2$: the generalisation of the proofs by M\"uller and Bajo is straightforward. \end{proof*}

\begin{remark*} The first proof for the invariance of the solvable radical is algebraic in nature and the second proof is geometric. \end{remark*} \newline

Finally, we consider the special case of nilpotent Lie algebras. In such a nilpotent Lie algebra $\Lg$, there are only two invariant subspaces: $0$ and $\Lg$ itself (proposition 3). From this it follows that the ideals $\Lg^k$, $\Lg^{(k)}$ and $\Lz_k(\Lg)$ of $\Lg$ are, in general, \emph{not} invariant. But the corresponding ideals (of the descending, derived or ascending series of the Leibniz-algebra) will be.

\section{Characterisation of nilpotent Lie algebras} We are now in a position to prove the main theorem. We include the proofs for completeness, even though their techniques are well-known. First, the general case is reduced to the solvable case. Then, the solvable case is reduced to the nilpotent case.

\begin{proposition*} If a Lie algebra has an invertible Leibniz-derivation, then the Lie algebra is solvable. \end{proposition*}

\begin{bewijs} Let $\Lg$ be a Lie algebra with an invertible Leibniz-derivation $P$, say of order $k$. We must then show that $\Lg$ is solvable. Consider a \emph{Levi-Mal'cev}-decomposition $\Lg = \Ls \ltimes \Lr$, where $\Lr$ is the solvable radical and $\Ls$ is a Levi-complement. Denote by $\map{\pi_{\Ls}}{\Lg}{\Ls}$ the projection along the solvable radical onto the Levi-complement. Then the endomorphism $Q = (\pi_{\Ls} \circ P) |_{\Ls}$ is a Leibniz-derivation of $\Ls$ of order $k$. Since $P$ is invertible and since $\Lr$ is invariant under $P$ (proposition $9$), we obtain $$ \Ls = \pi_{\Ls}(\Lg) = \pi_{\Ls} (P(\Lg)) = \pi_{\Ls} ( P (\Ls) ) + \pi_{\Ls} ( P (\Lr) ) = Q(\Ls).$$ This means that the Leibniz-derivation $Q$ is invertible. We may then invoke Jacobson's theorem or corollary $1$ to conclude that $\Ls = 0$. In other words: $\Lg$ is solvable.\end{bewijs}

\begin{lemma*}[Eigenspaces] Consider a Lie algebra $\Lg$ be with a Leibniz-derivation $P$ of order $k \geq 1$. Consider the generalised eigenspace decomposition of $\Lg$ with respect to the operator $P$: $ \Lg = \sum_{\alpha \in \Omega_P} \Lg_\alpha.$ Then for all generalised eigenvalues $\alpha_1,\ldots,\alpha_{k+1}$, we have the identity $$ \left[ \Lg_{\alpha_{1}},\ldots,\Lg_{\alpha_{k+1}} \right] \subseteq \Lg_{\alpha_{1} + \ldots + \alpha_{k+1}} .$$
\end{lemma*}

\begin{bewijs} For every $m \in \N$, $\alpha_1,\ldots,\alpha_{k+1} \in \bF$, and $x_1,\ldots,x_{k+1} \in \Lg$,
\begin{equation*}\left( P - \sum_{t = 0}^{k+1} \alpha_t \id_{\Lg} \right)^m ([ x_1,\ldots,x_{k+1} ])\end{equation*} is a linear combination of elements of the form
\begin{equation} [ (P - \alpha_1 \id_{\Lg})^{i_1}(x_1),\ldots, (P - \alpha_{k+1} \id_{\Lg})^{i_{k+1}}(x_{k+1})] \end{equation} for some $i_1,\ldots,i_{k+1}$ in $\N_0$ satisfying $i_1 + \ldots + i_{k+1} = m$. This can easily be shown using induction on $m$. Now assume that each $x_i$ is in $\Lg_{\alpha_i}$. It then suffices to find an $m$ for which all of the expressions $(1)$ vanish. For each generalised eigenvalue $\alpha_i$ there exists, by definition, an $m_i \in \N_0$ such that $(P - \alpha_i \id_{\Lg})^{m_i}(\Lg_{\alpha_i})$ vanishes. If we define $m$ to be $m_1 + \ldots + m_{k+1}$, then each element of the form $(1)$ will vanish, and this finishes the proof of the lemma. \end{bewijs}

\begin{proposition*} If a solvable Lie algebra has an invertible Leibniz-derivation, then the Lie algebra is nilpotent. \end{proposition*}

\begin{bewijs} We assume the notation of the previous proposition and lemma. Then for all $\alpha,\beta \in \Omega_P$ and $n \in \N_0$ we have $ \operatorname{ad}(\Lg_{\alpha})^{k n} (\Lg_{\beta}) \subseteq \Lg_{n k \alpha + \beta} $. We now choose an $n$, not depending on $\alpha$ or $\beta$, such that the right hand side of this inequality vanishes: $$ \operatorname{ad}(\Lg_{\alpha})^{k n} (\Lg_{\beta}) = 0. $$ This is possible since the set $\Omega_P$ of generalised eigenvalues is finite and since invertible transformations have only non-zero generalised eigenvalues. Any basis for $\Lg$ of generalised eigenvectors is then ad-nilpotent. Since $\Lg$ is solvable, \emph{Lie}'s theorem implies that \emph{all} elements of $\Lg$ are ad-nilpotent. \emph{Engel}'s theorem finishes the proof.\end{bewijs}

\begin{bewijs}[Main Theorem] Without loss of generality, we may assume that the base field is algebraically closed (cf. proposition $2.3$ in \cite{Burde}). The necessity is proved in Proposition 3. Together, propositions $11$ and $12$ prove the sufficiency. \end{bewijs}

We finish this paper with two remarks. \newline

\begin{remark*} Sufficient conditions on automorphisms and derivations for the nilpotency or solvability of a Lie algebra are also interesting and useful (but more complicated) in the case the base field has \emph{prime} characteristic. As an application, we cite the works by Shalev and Zel'manov on (pro) $p$-groups  of finite coclass, \cite{Shalev,ShalevZelmanov}. Note, however, that our main theorem can only hold in characteristic zero for the following trivial reason: \emph{every} Lie algebra over a field of characteristic $p > 0$ has an invertible $p$-derivation. The identity transformation is such a derivation. \end{remark*} \\

\begin{remark*} The so-called $n$-ary Lie algebras are natural generalisations of Lie algebras. One could suspect that a theorem, very much like Jacobson's, also holds for these generalised objects. It is important to note however, that this is \emph{not} the case: for each $n \geq 3$, there exist non-nilpotent $n$-ary Lie algebras with an invertible derivation. This has been shown in \cite{Williams}.  \end{remark*} \newline

\end{document}